\documentclass[11pt]{article}

\newtheorem{prop}{Proposition}[section]
\newtheorem{thm}{Theorem}[section]

\newtheorem{rem}{Remark}[section]
\newtheorem{defi}{Definition}[section]

\let\ch=\widehat
\let\fd=\rightarrow
\let\lfc=\mapsto
\let\lfd=\longrightarrow
\def\f#1{{\bf F}_{#1}}
\def\tr{\mathop{\rm Tr}\nolimits}
\let\dps=\displaystyle
\def\biq#1{\quad \hbox{ #1 }\quad }
\def\pt#1{\left(#1\right)}

\let\ss=\smallskip
\let\ms=\medskip

\let\bs=\bigskip
\let\bb=\bigbreak
\def\ssi{\hbox{ if and only if }}

\begin{document}

\title{Nonlinarity of Boolean functions\\ and  hyperelliptic curves}
\author{Eric F\'erard\thanks{Universit\'e de Polyn\'esie fran\c caise, Tahiti; {e-mail} {\tt ferard@upf.pf}}, 
Fran\c cois Rodier\thanks{Institut de Math\'ematiques de Luminy --
C.N.R.S. 
163 avenue de Luminy,
Case 907, Marseille Cedex 9, France; {e-mail} {\tt rodier@iml.univ-mrs.fr}}}
\date{}
\maketitle

Boolean functions is an important tool in computer sciences. It is especially useful in private key cryptography for designing stream  ciphers.
For security reasons, and also because Boolean functions
need also to have other properties than nonlinearity such as  balancedness or high
algebraic degree, it is important to have
 the possibility of choosing among many Boolean functions, not
only bent functions,
that is functions with the highest possible non linearity,
 but also functions which are  close to be bent in the sense that their nonlinearity is close to the nonlinearity of bent functions.
For $m$ odd, it would be particularly interesting to find functions with nonlinearity larger than the one of quadratic Boolean functions (called {\sl almost optimal} in \cite{cccf}). This has been done for instance
in the work of Patterson and  Wiedemann \cite{pw} and also of Langevin-Zanotti \cite{lz}.

\bb
Let $q=2^m$ and $\f{2^m}$ assimilated as a vector space to $\f2^m$.
In this talk, we want to study functions of the form $\tr G(x)$, where $G$ is a polynomial on $\f{2^m}$ and $\tr$ the trace of $\f{2^m}$ over $\f2$.

For $m$ even, many people got interested in finding bent functions of this form. To only mention the case  of monomials, one can get the known cases (Gold , Dillon/Dobbertin, Niho exponents) in the paper of  Leander \cite{le}. 

For $m$ odd, 
one might have expected that   among the functions 
$f:x\lfd \tr G(x)$
where $G$ is a polynomial of degree 7, there are some functions which are close to being bent in the previous sense.
This happens not to be the case, but
we will show that  for $m$ odd such functions have rather good nonlinearity or autocorrelation properties.
We use for that recent results of Maisner and Nart \cite{mn} about zeta functions of supersingular curves
of genus 2.

\bb
On the other hand,  vectorial Boolean functions are used in cryptography to construct block ciphers.
An
important criterion on these functions is a high resistance to the differential 
cryptanalysis. 
Nyberg \cite{ny} has introduced the notion of almost perfect nonlinearity (APN) to study differential attacks. 
We relate this notion to  the notion above, and
we will give some criterion for a function not to be almost perfect nonlinear.

\section{Preliminaries}

\subsection{Boolean functions}

Let $m$ be a positive integer and $q=2^m$.

\begin{defi}
A {Boolean function} with $m$ variables is a map
from the space
$V_m=\f2^m$ into $\f2$.
\end{defi}

A Boolean function is  linear if it is a linear form
on the vector space  $V_m$. 
It is  affine if it is
equal to a linear function up to addition of a constant. 


\subsection{Nonlinearity}
\label{defnl}

\begin{defi}
We call   nonlinearity of a Boolean function  
$f:V_m\lfd\f2$ 
the distance from $f$ to the set of affine functions  with
$m$ variables:
$$nl(f) = \min_{h \hbox{\,\scriptsize affine }} d(f, h)$$
where $d$ is the Hamming distance.
\end{defi}

One can show that the 
nonlinearity is equal to
$$\dps nl(f)  = 2^{m-1} - {1\over 2}\|\ch f\|_\infty
$$
where
$$  \|\ch f\|_\infty = 
\sup_{v\in V_m} \Bigl| \sum_{x\in V_m}\chi\pt{f(x)+v\cdot x}\Bigr|,$$
where $v\cdot x$ denote the usual scalar product   in $V_m$ and $\chi(f)=(-1)^f$.
It is the maximum of the Fourier transform of  $\chi(f)$ (the Walsh transform of $f$):
$$\ch f(v) = \sum_{x\in V_m}\chi\pt{f(x)+v\cdot x}.$$
Parseval identity can be written
$$\|\ch f\|_2^2= \frac{1}{q} \sum_{v \in V_m} \widehat{f}(v)^2=q$$
and we get,
for $f$ a Boolean function on $V_m$:
$$\sqrt q\le  \|\ch f\|_\infty\le q.$$

\subsection{The sum-of-square indicator}

Let $f$ be a Boolean function on $V_m$. Zhang and
 Zheng introduced the {\sl sum-of-square indicator}   \cite{zz}, as a measure of the {\sl global avalanche criterion}:
$$ \sigma_f={1\over q}{\sum_{x\in
V_m}\ch{f}(x)^4}=\|\ch{f}\|_4^4.$$
We remark that
\begin{equation}
\label{ineg2}
\|\ch f\|_2\le\|\ch f\|_4\le\|\ch f\|_\infty.
\end{equation}
Hence the values of $\|\ch f\|_4$ may be considered as a first approximation of  $\|\ch f\|_\infty$ and in some cases they may be easier to compute.
The relationship of this function with non-linearity was studied by A.~Canteaut et
al.\cite{cccf}.

\section{{The functions $f:x\lfd \tr\pt{ G(x)}$ where $G$ is a polynomial}}

\subsection{Divisibility of $\|\ch f\|_\infty$}

Let $G(x)$ be the polynomial $\sum_{i=0}^s a_ix^i$  with coefficients in $\f q$  and $f$ the Boolean function $\tr\circ G$. 

\begin{defi}
The binary degree of $G$ is the maximum value of $\sigma(i)$ for $0\le i\le s$, where $\sigma(i)$ is the sum of the  binary digits of $i$.
\end{defi}

One has the following proposition, due to C. Moreno and O. Moreno \cite{mm}.

\begin{prop}
Let $G$ be a polynomial with coefficients in $\f q$ and binary degree $d$.
Then $\|\ch f \|_\infty$ is divisible by $2^{\lceil{m\over d}\rceil}$.
\end{prop}

\subsection{Case where $G$ is a polynomial of binary degree 2}

The $\|\ch f\|_\infty$ are multiple of $2^{\lceil{m\over 2}\rceil}$. Therefore, if $m$ is even $\|\ch f\|_\infty$ is a multiple of $q^{1/2}$, and if $m$ is odd, of $\sqrt{2q}$. 
In particular, if $m$ is odd, the spectral amplitude is higher or equal to $\sqrt{2q}$ which is equal to that of the quadratic Boolean functions, of the maximum rank.

\section{The functions $f:x\lfd \tr\pt{G(x)}$ where $G$ is a binary polynomial of degree 3} 

One simply will study the case where $G$ is a binary polynomial of degree 2 to which one adds a monomial of degree 7: 
$$G=a_7x^7+\sum^s b_ix^{2^i+1}$$ 
where $a_7\ne0$ a polynomial of degree 7 with coefficients in $k$. We would like to evaluate $\|\ch f\|_4$ on $\f{2^m}$, for $f(x)=\tr\pt{G(x)}$ where $\tr$ indicates the function trace of $\f q$ on $\f2$: 
$$\tr(x)=\sum_{i=0}^{m-1}x^{2^i}.$$ 
One obtains the simple expression of $\|\ch f\|_4$ (cf \cite{ro1, ro2}):

\label{som}
$$\|\ch f\|_4^4 = \sum_{x_1+x_2+x_3+x_4=0}
\chi\pt{f(x_1)+f(x_2)+f(x_3)+f(x_4)}=q^2+\sum_{\alpha\in k^*}X_\alpha$$
with
$$ X_\alpha=\Big(\sum_{x\in k}\chi\circ\tr\pt{ G(x)+ G(x+\alpha)}\Big)^2.$$
To compute  $X_\alpha$, one can remark that
the curve of equation $y^2+y=G(x+\alpha)+G(x)$ is isomorphic to
\begin{eqnarray*}
\lefteqn{y^2+y =
G(\alpha)+ }\\
    &&+ (a_7 \alpha^6  
     + a_7^{1/4} \alpha^{3/4}  + a_7^{1/2}  \alpha^{5/2}
     +\sum (b_i\alpha)^{2^{-i}} +\sum b_i\alpha^{2^i} ) x +\cr 
 &&\quad  
   +( a_7 \alpha^4+a_7^{1/2} \alpha^{1/2}) x^3 + a_7 \alpha^2 x^5 
\end{eqnarray*}
which is an equation of a curve $C_1$ of genus  2 for $\alpha\ne0$.
One has
$$X_\alpha=(\#C_1-q-1)^2.$$

 To compute $X_\alpha$, we will need results of Van der Geer - van der Vlugt and  of Maisner - Nart.

\subsection{Van der Geer and van der Vlugt theory}

Let $C_1$ the curve with affine equation:
$$C_1 : y^2 + y = ax^5 + bx^3 + cx + d$$
with $a\ne0$.
Let $R$ be the linearized polynomial $ax^4 + bx^2 + c^2x$.
The map
\begin{eqnarray*}
Q : k &\fd& \f2\\ 
x &\lfc&  \tr(xR(x))
\end{eqnarray*}
 is the quadratic form associated to the  symplectic form
\begin{eqnarray*}
k \times k &\lfd& \f2\\ 
(x, y) &\lfc& < x, y >_R= \tr(xR(y) + yR(x)).
\end{eqnarray*}
The number of zeros of $Q$ determines the number of points of $C_1$:
$$\#C_1(k)=1+2\#Q^{-1}(0).$$
Let $W$ be the radical of the symplectic form $<,>_R$, and $w$ be its dimension over $\f2$. 
The codimension of the kernel $V$ of $Q$ in $W$ is equal to 0 or 1.

\begin{thm}
(van der Geer - van der Vlugt \cite{gv1}) 
\label{gv}

If $V \ne W$, then
$\#C_1(k) = 1 + q.$

If $V = W$, then
$\#C_1(k) = 1 + q \pm \sqrt{2^wq}.$

\end{thm}

    \subsection{Values of $X_\alpha$}

In \cite{fr}, we study the factorization of $P$ which determines $V$ and $W$ (see Maisner-Nart \cite{mn}).
   Thanks to the work of van der Geer - van der Vlugt, we can compute the number of points  of the curves $y^2+y=G(x+\alpha)+G(x)$.

 \begin{prop}
 \label{Xalpha}
Suppose that $m$ is odd. Then
$$
 X_\alpha=0 \biq{or} 2q  \biq{or} 8q.
$$
\ss
Let
$\ell=a_7 ^{-1/3}\alpha^{-7/3}$.
Then
$$\displaylines{
 X_\alpha=8q\ssi\hfill\cr
\tr\ell=0 \biq{,} \ell=v+v^4  \biq{with} \tr v=0\biq{,} 
\cr
\hfill
\tr\pt{{(a+c )\alpha\over\lambda} v^3}=1 \biq{,}
\tr\pt{{(a+c )\alpha\over\lambda} (v+v^2)}=1\quad;\cr
}$$
$$\displaylines{
 X_\alpha=2q\ssi
\tr\ell=1 \biq{;}
  \hfill
}$$
$$\displaylines{
 X_\alpha=0
\quad\hbox{in the remaining cases.}\hfill}$$
\end{prop}

\section{Evaluation of $\|\ch f\|_4^4$}

\begin{prop}
\label{eval}
The value of $\|\ch f\|_4^4$ on $\f{2^m}$ when $m$ is odd and $f(x)=\tr\pt{ G(x)}$ is such that
$$
|\|\ch f\|_4^4
-3q^2|\le 185.2^{s-1}q^{3/2}.$$
\end{prop}

Proof

One can evaluate the number of $\alpha$ which gives each case of the preceding proposition.
The proves of these evaluations are linked with the computations of exponential sums over the curve
$v+v^4=\gamma x^7.$
We get 
\begin{eqnarray*}
\Big |\#\{\alpha\mid X_\alpha=8q\}-{1\over 8}\Big|&\le&23.2^{s-1} q^{1/2}\cr
\noalign{\bs}
\Big |\#\{\alpha\mid X_\alpha=2q\}-{1\over 2}\Big|&\le&3q^{1/2}+1\cr
\end{eqnarray*}
One deduce easily the evaluation of $\|\ch f\|_4^4$.
The details of the proof will appear in \cite{fr}.

\begin{rem}
This result is to be compared with proposition 5.6 in \cite{ro1} where  the distribution of $\|\ch f\|_4^4$ for all Boolean function is shown to be concentrated around $3q^2$.
\end{rem}

\section{Bound for  $\|\ch f\|_\infty$}

From the theorem, we can deduce some lower bounds for $ \|\ch f\|_\infty$.
\begin{prop}
\label{borne}
For the functions $f:x\lfd \tr\pt{ G(x)}$ on $\f{2^m}$ where $G$ is the polynomial 
$G=a_7x^7+\sum^s b_ix^{2^i+1}$
 and $m$ is odd one has, for $m\le11+2s$:
$$\sqrt{2 q} \le \|\ch f\|_\infty.$$
For $m\ge15+2s$, one has moreover:
$$\sqrt{2 q}< \|\ch f\|_\infty.$$
\end{prop}

Proof

The evaluation of the number of $\alpha$ such that $\tr \ell=1$ in proposition \ref{Xalpha} gives:
$$2q^2-6q^{3/2}\le \|\ch f\|_4^4.$$
As it is easy to show that
$$\|\ch f\|_4^4\le q\|\ch f\|_\infty^2$$
we get
$2q-6q^{1/2}\le \|\ch f\|_\infty^2$
whence the result, as $\|\ch f\|_\infty$ is divisible by $2^{\lceil m/3\rceil}$.

The second inequality is a consequence of theorem \ref{eval}.

\begin{rem}
So $f$ is not almost optimal (in the sense of \cite{cccf}), for $m\ge15+2s$.
\end{rem}

\section{APN Functions}

Let us consider a function $G:\f q\lfd\f q$.

\begin{defi}
The function $G$ is said to be  {APN} (almost perfect nonlinear)
if for every $a\in\f q^*$ and $b\in \f q$, 
there exists at most 2 elements of $\f q$ such that\quad
$G(z+a)+G(z)=b$.
\end{defi}

\begin{prop}
The function 
$$\begin{array}{cccl}
G:  &\f q&\lfd&\f q      \\
\noalign{\ms}
  & x  &\lfc&  a_7x^7+\dps\sum_0^s b_ix^{2^i+1} 
  \end{array}
$$
  {is not APN} for $m\ge13+2s$.
\end{prop}

Proof

For $\gamma\in \f q$, consider the function\qquad
$f_\gamma(x)=\tr(G(\gamma x))$. The proposition follows from proposition \ref{eval} and the following result from Chabaud-Vaudenay \cite{cv}.
\begin{prop}
One has\quad
$\dps\sum_{\gamma\in k^*} \sigma(f_\gamma)\ge2q^2(q-1)$.

The function $G$ is  {APN} \ssi the  {\bf equality is true.}
\end{prop}

For $s\le2$, one can even say more.
The following theorem \cite{ro3} proves that the function $G$ is not APN for $m\ge 11$.

\begin{thm}
Let $G$ be a polynomial from $\f{2^m}$ to $\f{2^m}$, $d$ its degree. Let us suppose that the curve $X_\infty$ of equation
$${x_0^d+x_1^d+x_2^d+(x_0+x_1+x_2)^d\over (x_0+x_1)(x_2+x_1)(x_0+x_2)}=0$$
is smooth.
Then if
 {$m\ge6$ }
and  {$d<q^{1/6} +3.9$, }
 $G$ is not APN.
\end{thm}


\end{document}